\documentclass[a4paper,12pt]{article}
\usepackage{amsthm}
\usepackage{amssymb}
\usepackage{amsmath}
\usepackage{amsfonts}
\usepackage{color}
\usepackage{graphicx}

\newtheorem{remark}{Remark}[section]
\newtheorem{example}{Example}[section]

\newtheorem{lemma}{Lemma}[section]
\newtheorem{theorem}{Theorem}[section]
\newtheorem{proposition}{Proposition}[section]
\newtheorem{corollary}{Corollary}[section]

\def\b1{\mbox{\boldmath $1$}}

\newenvironment{demo*}{\vspace{3mm}\noindent{\bf Proof.}}{\hfill $\Box$ \vspace{3mm}}

\begin{document}
\title{\bf \Large    Joint mixability of elliptical  distributions and related families}
{\color{red}{\author{\normalsize{Chuancun Yin}
\\
{\normalsize\it  School of Statistics  and Data Science, Qufu Normal University}
\\\noindent{\normalsize\it Shandong 273165, China}\\e-mail:  ccyin@qfnu.edu.cn\\\and
\normalsize{Dan Zhu} 
\\
{\normalsize\it  School of Statistics  and Data Science, Qufu Normal University}\\
\noindent{\normalsize\it  Shandong 273165, China}\\
 }}}
\maketitle
\vskip0.01cm
\noindent{\large {\bf Abstract}}  { {
  In this paper,
three different proofs to a result of  Wang, Peng and   Yang (2013) which related to the joint mixability of  elliptical distributions with the same characteristic generator are present. Moreover, we generalize this result to any distributions with finite second moments. An open problem proposed by  Wang (2015) is solved by constructing  a   bimodal-symmetric distribution.
The joint mixability  of slash-elliptical distributions and skew-elliptical distributions is studied and the  extension to multivariate distributions   is also investigated.}}

\medskip





\numberwithin{equation}{section}
\section{Introduction}\label{intro}
A standard problem in variance reduction and simulation  is to minimize the variance of the sum $\sum_{i=1}^n X_i$ of random variables $X_1,\cdots, X_n$ with given marginal distributions $F_i$'s:
$$\min_{X_i\sim F_i}{\rm Var}\left(\sum_{i=1}^n X_i\right).$$
See Fishman (1972),  Gaffke and  R\"uscherndorf (1981) and Knott and Smith (2006) for more details. This problem has been considered, in a slightly different form, by R\"uschendorf and Uckelmann (1997, 2002). More generally, one solves the convex minimization
problem
$$\min_{X_i\sim F_i}{\rm E}f\left(\sum_{i=1}^n X_i\right),$$
 where $f$ is a convex function.  It is obvious that if there exist a constant $c$ such that $P(\sum_{i=1}^n X_i=c)=1$, then the  above two problems have optimal solutions. These solutions are related to the   concept of complete mixabilitiy   and joint mixability of distributions. For more details on the solutions of these problems, the relevant applications and a brief history of the concept of the
complete mixability  and joint mixability we refer the reader to  Wang and Wang (2011, 2016) and Wang et al. (2013), Puccetti and Wang (2015).
In recent years, the problem of studying complete mixability and joint mixability of distributions has received considerable attention. Complete mixability and joint mixability describe whether it is possible to generate random variables (or vectors) from given distributions with constant  sum.  The formally definition of complete mixability for a distribution   was  first introduced
in Wang and Wang (2011) and then extended to an
arbitrary set of distributions in Wang, Peng and Yang (2013), although
the concept has been used in variance reduction problems earlier (see Gaffke
and R\"uschendorf (1981), Knott and Smith (2006), R\"uschendorf and Uckelmann (2002)).
The properties are particularly of interest in quantitative
risk management  and optimization problems in the theory of optimal couplings, where dependence between risks is usually unknown or partially unknown.

Throughout the paper, we write $X\stackrel{d}{=}Y$ if the random variables (or vectors) $X$ and $Y$ have the same distribution. For a cumulative distributions function $F$, we write $X\sim   F$ to denote $F (x)=P(X\le x).$  By
convention, all vectors will be written in bold and will be considered as column vectors,
with the superscript $\top$ for transposition.  Next we introduce the concepts of completely mixable and jointly mixable distributions.
 Suppose $n$ is a positive integer.
We say the univariate distribution functions $F_1,\cdots, F_n$ are jointly mixable (JM) if there
exist $n$ random variables $X_1,\cdots,  X_n$ with distribution functions  $F_1,\cdots, F_n$, respectively,
such that
\begin{equation}
P(X_1 + \cdots + X_n =C) = 1,
\end{equation}
for some $C\in \Bbb{R}$. If (1.1) holds with $F_j = F, 1\le j\le n$, the distribution $F$ is said to be $n$-completely mixable ($n$-CM).  Any such $C$ is called a joint center of   $(F_1,\cdots, F_n)$. Accordingly, we say  $X_1,\cdots,  X_n$ are jointly mixed; see Wang, Peng and Yang (2013). Clearly, Equation (1.1) is equivalent to
$Var(X_1+\cdots+X_n)=0$ given the variance exists.
For  a brief history of the concept of the
complete mixability, we refer to   Wang  (2015) and Wang and Wang (2016). Existing results on complete
mixability and  joint mixability are summarized in Wang and Wang (2011),   Puccetti, Wang and Wang (2012),
 Wang and Wang (2016) and Puccetti et al. (2019).
As pointed out in Puccetti and  Wang (2015b),  as a full characterization of completely mixable distribution is still out of reach, there are even less results concerning
sufficient conditions for joint mixable distributions. The only available ones are given in the  paper of Wang and Wang (2016). Most studies in the above literature concerns
the shape of marginal density to justify the existence of joint mixability. Bignozzi and Puccetti (2015) extended the concept of joint mixability  and introduce the concept of $\phi$-joint mixability.
 Recent paper of Xiao and Yao (2020)  studies the dependence of joint mix random vectors from the perspective of covariance matrix. In this paper, we further develop the theory
of complete mixability and joint mixability for symmetric   distributions, elliptical distributions  and related families.

The rest of the paper is organized as follows. In Section 2 we discuss the conditions on a result of R\"uschendorf and Uckelmann (2002) related to complete mixability of
 continuous distribution function having a symmetric and unimodal density.
 Section 3 is dedicated to joint  mixability  of elliptical distributions and    slash/skew-elliptical distributions, respectively.
 Section 4 extended the result to the class of  multivariate  elliptical  distributions. Finally,  Section 5 gives a conclusion.

\vskip 0.1cm
\numberwithin{equation}{section}
 \section{Symmetric Distributions }
\setcounter{equation}{0}

It would be of interest to characterize the class of completely mixable distributions. Only partial characterizations are
known in the literature. One nice result for the complete mixability is given by R\"uschendorf and Uckelmann (2002), which is equivalent to state that any continuous distribution function having a symmetric and unimodal density is $n$-CM for any $n\ge 2$. Wang (2014) provided a new proof using duality representation. The property was also extended to multivariate distributions by R\"uschendorf and Uckelmann (2002).
\begin{lemma}\ (R\"uschendorf and Uckelmann (2002)). Any continuous distribution function
having a symmetric and unimodal density is $n$-CM for any $n\ge 2$.
\end{lemma}
We remark that the  inverse of Lemma 2.1  is not necessarily true. For example,  the density of Pearson type II distribution
 \begin{eqnarray}
f(x)=\left\{  \begin{array}{ll}\frac{1}{\pi\sqrt{1-x^2}},  \ &{\rm if}\ x\in (-1,1),\\
0, \ &{\rm if}\ x\notin (-1,1),\\
 \end{array}
  \right.
\end{eqnarray}
  is convex, bimodal and symmetric, so that
Lemma 2.1 can not applicable. Note that $f$ is
 $n$-CM for any  integer $n\ge 2$;  see  Puccetti, Wang and Wang (2012) for more details.
  Wang and Wang (2016) generalizes Lemma 2.1 and  studied  the joint mixability of unimodal-symmetric distribution based on a different technical
approach.
\begin{lemma} (Wang and Wang (2016)).
 Suppose that $F_1,\cdots,F_n$ are   distributions with unimodal-symmetric densities  from the same location-scale
family  with scale parameters $\theta_1,\cdots,\theta_n$, respectively. Then $F_1,\cdots,F_n$ is JM if and only if the scale inequality
\begin{equation}
\sum_{i=1}^n\theta_i\ge 2\max_{1\le i\le n}\theta_i
\end{equation}
 is satisfied.
\end{lemma}
{\bf Proof} \ Using Theorem 3.1 in Wang and Wang (2016), we can give a more simple proof to ``if part". In fact, $X_i\sim F_i$ can be written as $X_i\stackrel{d}{=} \theta_i RU_i+\mu$, where
$\mu$ is a constant, $R$ is a random variable on $(-\infty,\infty)$ and  $U_i$ is uniformly distributed on $(-1,1)$ independent of $R$. The result follows from  Theorem 3.1 in Wang and Wang (2016)  since  $\theta_iU_i\sim U(-\theta_i,\theta_i)$, $i=1, \cdots,n$.  This ends the proof.

Suppose that $Y$ has a  distribution function $F$   and that $\theta$ has a distribution function $H$  on $(0, \infty)$ and, $Y$ and $\theta$ are independent. Then the distribution of $X = \theta Y$ is referred
to as a scale mixture of $F$ with a scale mixing distribution $H$.
The following corollary  is a direct consequence of Lemma 2.1.
\begin{corollary}
The scale mixture of  an unimodal-symmetric continuous distribution  with center $\mu$ is $n$-CM ($n\ge 2$) with center $\mu$.
\end{corollary}

 The complete mixability and joint mixability   is a concept of negative dependence (cf. Puccetti and Wang (2015a)) and not all univariate distributions $F$ are $n$-CM.
 If the supports of $F_i$ $(i=1,2,\cdots, n$) are unbounded from one side,   then  ($F_1,\cdots, F_n$)   is not  JM for any $n\ge 1$;   see Remark 2.2 in Wang and Wang (2016).
 Now we list more examples (The proof learned largely from Ruodu Wang).

 \begin{example}\ Assume $F_1, \cdots, F_{2n+1}$ are  $2n+1$  univariate  distribution functions with symmetric densities  on the same interval $[-a,a]$ $(a>0)$, if $F_i(\frac{n}{n+1}a)\le \frac{n+1}{2n+1}$ $(i=1,2,\cdots, 2n+1)$, then  $(F_1,\cdots, F_{2n+1})$   is not  JM.
\end{example}

{\bf Proof}\ For any $X_i\sim F_i$ $(i=1,2,\cdots,2n+1)$, the conditions $F_i(\frac{n}{n+1}a)\le \frac{n+1}{2n+1}$ $(i=1,2,\cdots, 2n+1)$ imply that
$$P\left(|X_i|>\frac{n}{n+1}a\right)>\frac{2n}{2n+1},\; i=1,2,\cdots, 2n+1.$$
It follows that
\begin{eqnarray*}
P\left(|X_1|> \frac{na}{n+1},\cdots, |X_{2n+1}|> \frac{na}{n+1}\right)&\ge& \sum_{i=1}^{2n+1}P\left(|X_i|>  \frac{na}{n+1}\right)-2n\\
&>& (2n+1)\frac{2n}{2n+1}-2n=0.
\end{eqnarray*}
Note that
$$\left\{\sum_{i=1}^{2n+1}X_i\neq 0\right\}\supseteq \left\{|X_1|> \frac{na}{n+1},\cdots, |X_{2n+1}|> \frac{na}{n+1}\right\}.$$
Hence
$$P\left(\sum_{i=1}^{2n+1}X_i\neq 0\right)>0.$$
Thus ($F_1,\cdots, F_{2n+1}$)   is not  JM.

\begin{corollary}(Necessary Condition)\ Assume $F_1, \cdots, F_{2n+1}$ are  $2n+1$  univariate  distribution functions with symmetric densities  on the same interval $[-a,a]$ $(a>0)$, if $(F_1,\cdots, F_{2n+1})$   is   JM, then there exists some $i$ $(1\le i\le 2n+1)$ such that
$$P\left(|X_i|\le \frac{n}{n+1}a\right)>\frac{1}{2n+1}.$$
\end{corollary}

 The interval $[-a,a]$ in Example 2.1 can be changed as $(-\infty,\infty)$.

 \begin{example}\ Assume $F_1, \cdots, F_{2n+1}$ are  $2n+1$  univariate  distribution functions with symmetric densities  on the same interval $(-\infty, \infty)$, if
  there exists $a>0$ such that
  $$F_i(a)-F_i\left(\frac{n}{n+1}a\right)\ge \frac{n}{2n+1},\;\; i=1,2,\cdots, 2n+1,$$
   then  $(F_1,\cdots, F_{2n+1})$   is not  JM.
\end{example}
{\bf Proof}\ For any $X_i\sim F_i$ $(i=1,2,\cdots,2n+1)$, the conditions
$$F_i(a)-F_i\left(\frac{n}{n+1}a\right)\ge \frac{n}{2n+1},\;\; i=1,2,\cdots, 2n+1,$$
 imply that
$$P\left(|X_i|\in \left[\frac{n}{n+1}a,a\right]\right)>\frac{2n}{2n+1},\; i=1,2,\cdots, 2n+1.$$
It follows that
\begin{eqnarray*}
P\left(\bigcap_{i=1}^{2n+1}\left\{|X_i|\in \left[\frac{na}{n+1},a\right]\right\}\right)&\ge& \sum_{i=1}^{2n+1}P\left(|X_i|\in \left[\frac{na}{n+1},a\right]\right)-2n\\
&>& (2n+1)\frac{2n}{2n+1}-2n=0.
\end{eqnarray*}
Note that
$$\left\{\sum_{i=1}^{2n+1}X_i\neq 0\right\}\supseteq \bigcap_{i=1}^{2n+1}\left\{|X_i|\in \left[\frac{na}{n+1},a\right]\right\}.$$
Hence
$$P\left(\sum_{i=1}^{2n+1}X_i\neq 0\right)>0.$$
Thus ($F_1,\cdots, F_{2n+1}$)   is not  JM.

 The following example  tells us the symmetry of $F$ does not implied $F$ is 3-CM and  that the   unimodality assumption on the density   can not removed.

{\bf Example 2.3} Assume  that $F$  has the following bimodal symmetric density
  \begin{eqnarray*}
f(x)=\left\{  \begin{array}{ll} \frac{2r+1}{2a^{2r+1}} x^{2r},  \ &{\rm if}\ x\in [-a,a],\\
0, \ &{\rm if}\ x\notin [-a,a],\\
 \end{array}
  \right. \nonumber
\end{eqnarray*}
where $r$ is a positive integer.
The distribution is given by
\begin{eqnarray*}
F(x)=\left\{  \begin{array}{lll} 0,  \ &{\rm if}\ x<-a,\\
\frac{1}{2a^{2r+1}}\left(x^{2r+1}+a^{2r+1}\right), \ &{\rm if}\ -a\le x<a,\\
1,\ &{\rm if}\ x\ge a.\\
 \end{array}
  \right. \nonumber
\end{eqnarray*}
It is easy to see that
 $$F\left(\frac{a}{2}\right)=\frac{1}{2}+\frac{1}{2^{2r+2}}<\frac{2}{3}.$$
Thus $(F,F,F)$  is not  JM.
Or, equivalently,   $F$ is not 3-CM.

{\bf Example 2.4} Assume  that $F$  has the following   bimodal symmetric density
  \begin{eqnarray*}
f_m(x)=\left\{  \begin{array}{ll} C_m\frac{x^{2m}}{\sqrt{1-x^2}},  \ &{\rm if}\ x\in (-1,1),\\
0, \ &{\rm if}\ x\notin (-1,1),\\
 \end{array}
  \right.
\end{eqnarray*}
where $C_m$ is a  normalizing constant and $m\ge 0$ is an integer.
When $m=0$, $f_m$ is the density of Pearson type II distribution, thus it is $n$-CM for any  integer $n\ge 2$.
When $m\ge 1$, using Example 2.1 one can check that  $F$ is not $(2n+1)$-CM for not very large $n$.
  The $(2n+1)$-complete mixability of  this distribution is not covered by any known theoretical results for large $n$.
 Further we  consider a distribution $F$ with density
  $$f(x)=\sum_{m=1}^{\infty}\alpha_m f_m(x),$$
   where $\{\alpha_m\}_{m\ge 1}$ is a sequence of positive values with $\sum_{m=1}^{\infty}\alpha_m=1$.
   This distribution is not $(2n+1)$-completely mixable  for any $n\ge 1$.
   Thus we give a counterexample   to the following open problem:

 {\bf Open Problem}  (Wang (2015)).  Are all absolutely continuous distributions on a bounded interval $n$-CM
for large enough $n$?

 \vskip 0.1cm
 \numberwithin{equation}{section}
 \section{Elliptical Distributions  and Related Families}
\setcounter{equation}{0}

\subsection{Elliptical Distributions }

Let ${\bf \Psi}_n$ be a class of functions $\psi: [0,\infty) \rightarrow \mathbb{R}$
such that function $\psi(|\bf t|^2), t\in \Bbb{R}^n$  is an $n$-dimensional
characteristic function. It is clear that
${\bf \Psi}_n\subset {\bf \Psi}_{n-1}\cdots\subset {\bf \Psi}_1.$
Denote by ${\bf \Psi}_{\infty}$  the set of characteristic generators that generate an
$n$-dimensional elliptical distribution for arbitrary $n\ge 1$. That is
${\bf \Psi}_{\infty}=\cap_{n=1}^{\infty}{\bf \Psi}_{n}.$

{\bf Defination 3.1} A random vector ${\bf{X}}=(X_1,X_2,\cdots,X_n)^{\top}$ is said to have an elliptical distribution with
parameters $\boldsymbol{\mu}$  and $\bf{\Sigma}$, written as $\bf{X}\sim {\bf E}_n({\boldsymbol \mu},{\bf \Sigma}, {\bf \psi})$, if its characteristic function can be expressed as
\begin{equation}
{\bf \varphi_{\bf X}(t)}=\exp\left(i{\bf t}^{\top}{\boldsymbol{\mu}}\right)\psi\left({\bf t}^{\top}\bf{\Sigma}{\bf t}^{\top}\right), \;{\bf t}\in\Bbb{R}^n,
\end{equation}
for some column-vector $\boldsymbol{\mu}$, $n\times n$ positive semidefinite
matrix $\bf{\Sigma}$ and for some function $\psi\in {\bf \Psi}_n$ with $\psi(0)=1$, which
is called the characteristic generator. In general, elliptical distributions can be bounded or unbounded, unimodal or multimodal.
 When $\psi(u)=\exp(-u/2)$, ${\bf E}_n({\boldsymbol \mu},{\bf \Sigma}, {\bf \psi})$ is the normal distribution  ${\bf N}_n({\boldsymbol \mu},{\bf \Sigma})$  and when $n=1$ the class of elliptical distributions
consists  of the class of symmetric distributions.
 It is well known that an $n$-dimensional random vector $\bf{X}\sim {\bf E}_n({\boldsymbol \mu},{\bf \Sigma}, {\bf \psi})$
 if and only if for any vector ${\boldsymbol{\alpha}}\in \Bbb{R}^n$, one has (Cambanis et al. (1981))
 ${\bf\boldsymbol\alpha^{\top} X}\sim {\bf E}_1({\boldsymbol \alpha^{\top}\boldsymbol\mu}, \boldsymbol\alpha^{\top}{\bf \Sigma}\boldsymbol\alpha, {\bf \psi}).$
In particular, $X_i\sim {\bf E}_1(\mu_i, \sigma_i^2,  \psi)$ and $\sum_{i=1}^n X_i\sim {\bf E}_1({\boldsymbol e_n^{\top}\boldsymbol\mu}, \boldsymbol e_n^{\top}{\bf \Sigma}\boldsymbol e_n, {\bf \psi})$.

The next result is due to Wang, Peng and   Yang (2013), see also  Wang and Wang (2016). Here we present three another proofs by finding the exact  dependence structure and by using Lemma 2.2, respectively.
\begin{theorem}  (Wang, Peng and   Yang (2013))\;
Suppose  $F_i\sim   {\bf E}_1(\mu_i, \sigma^2_i,  \psi)$, where $\mu_i\in {\Bbb R}, \sigma_i>0,$ $\psi$  is a characteristic generator
for an $n$-elliptical distribution.   Then ($F_1,\cdots, F_n$)   is JM  if and only if
\begin{equation}
\sum_{i=1}^n \sigma_i\ge 2\max\{\sigma_1,\cdots,\sigma_n\}.
\end{equation}
\end{theorem}
{\bf Proof}\ The proof of the only-if-part is the same as that of Wang, Peng and   Yang (2013) with a minor revision.  Without loss of generality, we assume $\sigma_1\ge \sigma_2\ge\cdots\ge\sigma_n$. If there exist $X_1\sim F_1,\cdots, X_n\sim F_n$ such that Var($X_1+\cdots+X_n)=0$, then
\begin{eqnarray*}
Var(X_1+\cdots+X_n)&=&Var(X_1)+Var(X_2+\cdots+X_n)\\
&&+2Cov(X_1,X_2+\cdots+X_n)\\
&&\ge \left(\sqrt{Var(X_1)}-\sqrt{Var(X_2+\cdots+X_n)}\right)^2,
\end{eqnarray*}
which implies
\begin{eqnarray*}
0=\sqrt{Var(X_1+\cdots+X_n)}\ge \sqrt{Var(X_1)}-\sqrt{Var(X_2+\cdots+X_n)},
\end{eqnarray*}
from which we get
$$\sum_{i=2}^n \sigma_i\ge \sigma_1,$$
as desired.

For the if-part we will present three another proofs.

{\bf First proof}. Assume  ${\bf{X}}\sim {\bf E}_n({\boldsymbol \mu},{\bf \Sigma}, {\bf \psi})$,
where ${\bf\boldsymbol \mu}=(\mu_1,\cdots,\mu_n)^{\top}$ and ${\bf \Sigma}=(\sigma_{ij})_{n\times n}$.
Here
\begin{eqnarray*}
  \sigma_{ij}=\left\{  \begin{array}{ll} \sigma_i^2,  \ &{\rm if}\ i=j,\\
\frac{1}{(n-1)(n-2)}(\sigma^2_k-\sum_{l\neq k}\sigma^2_l), \ &{\rm if}\ k\neq i\neq j. \end{array}
  \right.
\end{eqnarray*}
It is straightforward to check that ${\bf \Sigma}$ is positive semidefinite under condition (3.2) and
the summation of all entries in   ${\bf \Sigma}$ is zero. Each component  $X_i$ of  $\bf X$  has distribution
${\bf E}_1(\mu_i, \sigma_i^2, \psi)$, $i=1,2,\cdots, n$.
The characteristic function of $\sum_{i=1}^n X_i$ can be expressed as
\begin{equation}
\varphi_{\sum_{i=1}^n X_i }(t)=\exp\left(it\sum_{i=1}^n\mu_i\right)\psi(0)=\exp\left(it\sum_{i=1}^n \mu_i\right), \;t\in\Bbb{R}.
\end{equation}
  Hence,
$$P\left(\sum_{i=1}^n X_i= \sum_{i=1}^n \mu_i\right)=1,$$ and thus ($F_1,\cdots, F_n$)   is JM.

{\bf Second proof}.  Considering  the  same  ${\bf \Sigma}=(\sigma_{ij})_{n\times n}$ as in the first proof, obviously, the  summation of   each row  in  ${\bf \Sigma}$  is zero. The  if-part follows from  Proposition 5 in Xiao and Yao (2020) which is says that if $(X_1,X_2,\cdots,X_n)'$ has a covariance matrix ${\bf \Sigma}$, then, it is a joint mix if and only if each row sum of ${\bf \Sigma}$ is 0.

{\bf Third proof}. We remark that if  $F_i\sim   {\bf E}_1(\mu_i, \sigma^2_i,  \psi)$ has a density and $\psi$  is a characteristic generator
for an $n$-elliptical distribution $(n\ge 2)$, then $F_i$  is unimodal and symmetric. Thus the  if-part   is a direct consequence of Lemma 2.2. $\hfill\square$

From the second proof of Theorem 3.1, it can be seen that the conclusion of the theorem also holds for general random variables. So we have the following proposition.
\begin{proposition}
Suppose distributions $F_i$ have finite variances $\sigma_i>0$, $i=1,2,\cdots,n$.    Then ($F_1,\cdots, F_n$)   is JM  if and only if
\begin{equation*}
\sum_{i=1}^n \sigma_i\ge 2\max\{\sigma_1,\cdots,\sigma_n\}.
\end{equation*}
\end{proposition}

The following result is a direct consequence of Theorem 3.1.
\begin{corollary}
Suppose  $F\sim   {\bf E}_1(\mu, \sigma^2,  \psi)$ with  $\psi\in {\bf \Psi}_{\infty}$. Then $F$ is $n$-CM for any $n\ge 2$.
\end{corollary}

 \begin{remark} Theorem 2.21 in Fang, Kotz and Ng (1990) shows that  $\psi\in {\bf \Psi}_{\infty}$ if and only if  $F\sim   {\bf E}_1(\mu, \sigma^2,  \psi)$  is a mixture of normal distributions. Some such elliptical distributions are  normal distribution,   $T$-distribution, Cauchy distribution, stable laws distribution and Pearson type VII distribution; see Andrews and Mallows (1974) and Kano (1994).
\end{remark}

 Note that there are  continuous, unimodal and symmetric densities do not belong to the class of normal scale mixtures; see West (1987). Thus Corollary 3.1 is a special case of Corollary 2.1. In the sequel, we list more examples.

  {\bf Example 3.1}\ Consider the generalized logistic distribution with density
  \begin{equation}
  f(x)=C\frac{\exp(-\alpha x^{\beta})}{(1+\exp(- x^{\beta}))^{2\alpha}}, -\infty<x<\infty,
  \end{equation}
  where $C>0, \alpha>0,\beta>0$ are constants. If $\alpha=1$ and $\beta=2$, (3.4) is 1-dimensional logistic distribution which  is  unimodal and symmetric but not a scale mixture of normal densities; see  G\'omez-S\'anchez-Manzano  et al. (2006).  If $\alpha=\beta=1$, (3.4) is  standard logistic distribution  which  is  unimodal and symmetric and  can be represented as  a scale mixture of normal densities; see Stefanski (1990).

 {\bf Example 3.2}\
 Kotz type distributions  with  density generator
$$g(r)=C r^{N-1}\exp(-mr^{\beta}),\ m, \beta>0, N>1$$
 have symmetric and bimodal  densities,
  the $(2n+1)$-complete mixability of  those  Kotz type distributions is not covered by any known theoretical results for  $n\ge 1$.

\numberwithin{equation}{section}
 \subsection{Slash-Elliptical Distributions }
In this subsection, we investigate  joint  mixability  of  slash-elliptical distributions.  We say that a random variable $X$ follows a slash elliptical  distribution if it can
be written as
\begin{equation}
X=\frac{Z}{U^{\frac{1}{q}}}+\mu,
\end{equation}
 where $Z\sim   {\bf E}_1(0, \sigma^2,  \psi)$  is independent of $U\sim U(0,1)$ and $q > 0$ is the parameter related to the distribution
kurtosis. We use the notation $X \sim  {\bf SE}_1(\mu, \sigma^2,  \psi;q)$.
  Similarly, we say that a random vector ${\bf X}\in {\Bbb R}^p$  has slash-elliptical
multivariate distribution with vector location parameter ${\boldsymbol\mu}$, positive semidefinite matrix scale parameter ${\bf \Sigma}$, and tail
parameter $q>0$, if it can be represented as
 \begin{equation}
 {\bf X}=\frac{\bf Z}{U^{\frac{1}{q}}}+{\boldsymbol\mu},
 \end{equation}
 where ${\bf Z}\sim   {\bf E}_p(0, {\boldsymbol\Sigma},  \psi)$  is independent of $U\sim U(0,1)$ and kurtosis
parameter  $q > 0$.  We denote this as
 ${\bf X} \sim {\bf SE}_p({\boldsymbol \mu},{\bf \Sigma}, {\bf \psi};q)$.
 Properties of this family are discussed in G\'omez, Quintana and Torres (2007) and  Bulut and Arslan (2015).

Using the representation (3.5), the following theorem is a consequence  of Theorem 3.1.
\begin{theorem}
Suppose that $F_i\sim  { \bf SE}_1(\mu_i,\sigma_i^2,  \psi;q)$, where $\mu_i\in {\Bbb R}, \sigma_i>0$,  $\psi$  is a characteristic generator
for an $n$-variate slash-elliptical
 distribution.   Then ($F_1,\cdots, F_n$)   is JM   if and only if
 $$\sum_{i=1}^n \sigma_i\ge 2\max\{\sigma_1,\cdots,\sigma_n\}.$$
\end{theorem}


 \subsection{Skew-Elliptical Distributions }

 A  univariate random variable ${Z}$ has a  skew-elliptical
 distribution  if its probability density function (pdf) is
$$2g(z)\pi(\lambda z), -\infty < z < \infty, $$
where  $g$  is a pdf of univariate elliptical distribution  with center $0$,  $\lambda \in  {\mathbb R}$ and $\pi$ is the distribution   function of $g$. We write $Z\thicksim SE_1(0, g, \pi,\lambda)$. In particular, if  $g$ is the pdf of $N(0,1)$, then
${Z}$ is called has a  skew-normal  distribution and write $Z\thicksim  SN(\mu, \sigma^2, \lambda)$.
A random variable $X$ follows a skew scale mixture of normal   distribution with location parameter $\mu \in \mathbb{R}$, scale parameter  $\sigma^2$ and skewness parameter $\lambda\in \Bbb{R}$  if its pdf is given by
$$g(x)=2g_0(x)\Phi \left(\lambda \frac{x-\mu}{\sigma}\right), \; x\in {\mathbb R},$$
where $\lambda \in {\Bbb R}$, $\Phi$ is   the distribution function of $N(0, 1)$ and  $g_0$ is the pdf of scale mixture of normal  distribution defined as
$$g_0(x)=\int_0^{\infty}\phi(x;\mu,v^2\sigma^2)dH(v).$$
Here $H$ is a (unidimensional) probability distribution function such that $H(0) = 0$. We use the notation $X\sim SSMN(\mu,\sigma^2, \lambda, H)$. For more details see Andrews and Mallows (1974) and G\'omez-S\'anchez-Manzano et al. (2006).

 For complete mixability of skew-normal distribution, we have the following result.
\begin{theorem}
Suppose that $F$ has the distribution $SN(\mu, \sigma^2, \lambda)$,    then $F$ is not $n$-CM ($n\ge 2$)  for sufficiently large  $|\lambda|$.
\end{theorem}
{\bf Proof} \ Assume random variable $X_{\lambda}$ has the distribution $F$.  Since  $X_{\lambda}  \stackrel{d}{=} \sigma X+\mu$, where $X\sim SN(0, 1, \lambda)$, we prove the theorem for the case  $X_{\lambda}\sim SN(0, 1, \lambda)$ only.
 It follows from Henze (1986) that
if  $X_{\lambda}\sim SN(0, 1, \lambda)$, then it has the stochastic representation
$$X_{\lambda}  \stackrel{d}{=} \frac{\lambda}{\sqrt{1+\lambda^2}}|U|+\frac{1}{\sqrt{1+\lambda^2}}V,$$
where $U$ and $V$ are independent $N(0, 1)$ random variables. Moreover,
$$E(X_{\lambda})=\frac{\lambda}{\sqrt{1+\lambda^2}}\sqrt{\frac2\pi}.$$
 Without loss of generality, we assume $\lambda>0$.
For any   $U_i$ and $V_i$ are independent $N(0, 1)$ random variables, we have
\begin{eqnarray*}
&P&\left(\frac{\lambda}{\sqrt{1+\lambda^2}}\sum_{i=1}^n |U_i|+\frac{1}{\sqrt{1+\lambda^2}} \sum_{i=1}^n V_i> nE(X_{\lambda})\right)\\
&&\ge P\left(\frac{\lambda}{\sqrt{1+\lambda^2}}\sum_{i=1}^n |U_i|+\frac{1}{\sqrt{1+\lambda^2}} \sum_{i=1}^n V_i > nE(X_{\lambda}), \bigcap_{i=2}^n\{\frac{\lambda}{\sqrt{1+\lambda^2}} |U_i|+\frac{1}{\sqrt{1+\lambda^2}} V_i\ge 0\} \right)\\
&&\ge P\left(\frac{\lambda}{\sqrt{1+\lambda^2}}|U_1|+\frac{1}{\sqrt{1+\lambda^2}}V_1 > nE(X_{\lambda}), \bigcap_{i=2}^n\{\frac{\lambda}{\sqrt{1+\lambda^2}} |U_i|+\frac{1}{\sqrt{1+\lambda^2}} V_i\ge 0\}\right)\\
&&\ge 1-P\left(\frac{\lambda}{\sqrt{1+\lambda^2}}|U_1|+\frac{1}{\sqrt{1+\lambda^2}}V_1\le nE(X_{\lambda})\right)\\
&&-\sum_{i=2}^n P\left(\frac{\lambda}{\sqrt{1+\lambda^2}} |U_i|+\frac{1}{\sqrt{1+\lambda^2}} V_i<0\right)>0,
\end{eqnarray*}
 for sufficiently large  $\lambda$.
This shows that the distribution $SN(0, 1, \lambda)$   is not $n$-CM ($n\ge 2$)  for sufficiently large  $|\lambda|$.  $\hfill\square$
\begin{remark}\  For $F\sim SN(\mu, \sigma^2, \lambda)$,  we conjecture that there exists an integer $n_0(\lambda)$ such that  $F$ is not $n$-CM for $n\le n_0(\lambda)$ and,  $F$ is $n$-CM for $n> n_0(\lambda)$; For an integer $n\ge 2$, there exists a $\lambda_0(n)\ge 0$ such that $F$ is $n$-CM whenever $|\lambda|\in [0, \lambda_0(n)]$ and, $F$ is not  $n$-CM whenever $|\lambda|\in (\lambda_0(n),\infty)$.
\end{remark}
\begin{theorem}
Suppose that $F$ has the distribution  $SSMN(\mu,\sigma^2, \lambda, H)$,  then $F$ is not $n$-CM ($n\ge 2$)  for sufficiently large  $|\lambda|$.
\end{theorem}
{\bf Proof} For any $X_i\sim SSMN(\mu,\sigma^2, \lambda, H)$ ($i=1,2,\cdots,n$) and $V\sim H$ such that $V$ is independent of $X_1,\cdots, X_n$.
By the definition of skew scale mixture of normal   distribution, we have
$$X_i|V=v\sim SN(\mu, \sigma^2 v^2, \lambda v).$$
Then for any constant $C$ and   sufficiently large  $|\lambda|$,
using Theorem 3.3,
 $$P\left(\sum_{i=1}^n X_i=C\right)=\int_0^{\infty}P\left(\sum_{i=1}^n X_i=C|V=v\right)dH(v)<\int_0^{\infty} dH(v)=1.$$
Thus  $F$ is not $n$-CM  for sufficiently large  $|\lambda|$. $\hfill\square$

\begin{remark} It seems  we can guess that as long as $F$ is asymmetric on $(-\infty,\infty)$ with unbounded support from two sides, then  $F$ is not  $n$-CM. But it is wrong.
  The following is a counterexample.  Assume   $P$  is continuous distribution on  interval $(-1, 1)$ having an asymmetric concave density and centered at 0,   $Q$ is normal $N(0,1)$.
Then for any $\lambda \in (0, 1), \lambda P + (1 - \lambda)Q$ is asymmetric  and  by the additivity (see Proposition 2.1 (3) in Wang and Wang (2011)) it  is $n$-CM for $n\ge 3$.
\end{remark}

\vskip 0.2cm
 \section{ Extensions to multivariate distributions}
\setcounter{equation}{0}

In this section we extent some results in Section 3  to the  class of  $n$-variate elliptically contoured distributions.
We first introduce some notions. The notation $vec(\bf A)$ denotes the vector $(\bf a_1^{\top},\cdots, a_n^{\top})^{\top}$, where ${\bf a_i}$ denotes the $i$th column of $p\times n$ matrix
${\bf A}$, $i=1,2,\cdots, n$. we use ${\bf A}\otimes {\bf B}$ to denote the  Kronecker product of the matrices ${\bf A}$ and ${\bf B}$; We use $tr({\bf A})$ to denote the  trace of the square matrix ${\bf A}$ and  $etr({\bf A})$ to denote  $exp(tr({\bf A}))$ if ${\bf A}$ is a square matrix.
We  use the following definition given in Gupta,  Varga and Bodnar (2013).

{\bf Definition 4.1}. Let ${\bf X}$ be a random matrix of dimensions $p\times n$. Then, ${\bf X}$ is said to
have a matrix variate elliptically contoured distribution if its characteristic
function has the form
\begin{equation}
E(etr(i{\bf T}^{\top}{\bf X}))=etr(i{\bf T}^{\top}{\bf M})\Psi(tr(\bf T^{\top}\Sigma T\Phi)).
\end{equation}
with ${\bf T}: p\times n, {\bf M}: p\times n, {\bf\Sigma}: p\times p, {\bf \Phi} : n\times n, {\bf \Sigma}\ge 0$ (positive semidefinite), ${\bf \Phi}\ge 0$, and $\Psi: [0,\infty)\rightarrow{\Bbb R}$.
This distribution will be denoted by ${\bf X}\sim E_{p,n}({\bf M}, {\bf \Sigma}\otimes {\bf \Phi}, \Psi)$.

The important special case of  matrix variate elliptically contoured distribution is the
matrix variate normal distribution (${\bf X}\sim N_{p,n}({\bf M}, {\bf \Sigma}\otimes {\bf \Phi})$), its characteristic function is
\begin{equation}
E(etr(i{\bf T}^{\top}{\bf X}))=etr\left(i{\bf T}^{\top}{\bf M}-\frac12\bf T^{\top}\Sigma T\Phi\right).
\end{equation}
The next lemma shows that linear functions of a random matrix with  matrix variate elliptically contoured
distribution have elliptically contoured distributions also (see Theorem 2.2 in Gupta,  Varga and Bodnar (2013)).
\begin{lemma}
Let ${\bf X}\sim E_{p,n}({\bf M}, {\bf \Sigma}\otimes {\bf \Phi}, \Psi)$.
 Assume ${\bf C} : q\times m, {\bf A} : q\times p$, and  ${\bf B} :
n\times m$ are constant matrices. Then,
$${\bf AXB+C} \sim E_{p,n}({\bf AMB+C}, {\bf A\Sigma A^{\top}}\otimes {\bf B^{\top}\Phi B}, \Psi).$$
\end{lemma}

The next lemma gives the marginal distributions of a  matrix variate elliptically contoured
distribution  (see Theorem 2.9 in Gupta,  Varga and Bodnar (2013)).
\begin{lemma}
Let ${\bf X}\sim E_{p,n}({\bf M}, {\bf \Sigma}\otimes {\bf \Phi}, \Psi)$, and partition ${\bf X, M}$, and ${\bf \Phi}$ as
$${\bf X}= ({\bf X_1, X_2}), {\bf M}=({\bf M_1}, {\bf M_2}),$$
and
$${\bf \Phi}=\left[\begin{array}{cc} \Phi_{11}&\Phi_{12}\\
\Phi_{21}&\Phi_{22}
\end{array}
\right],$$
where ${\bf X_1}$ is $p\times m, {\bf M_1}$ is $p\times m$, and ${\bf \Phi_{11}}$ is $m\times m, 1\le m < n.$ Then
 $${\bf X_1}\sim E_{p,m}({\bf M_1}, {\bf \Sigma}\otimes {\bf \Phi_{11}}, \Psi).$$
\end{lemma}

\begin{theorem}
Suppose that $F\sim   {\bf E}_p({\bf 0},  {\bf \Sigma},  \Psi)$, where ${\bf \Sigma}\ge 0 $ is a  $p\times p$ matrix, $\Psi$  is a characteristic generator
for a $p\times n$  matrix variate elliptically contoured distribution  ($n\ge 2$). Then  there exist
$n$ $p$-dimensional random vectors ${\bf X}_1,\cdots, {\bf X}_n$ identically distributed as $F$  such that
 $$P({\bf X}_1 + \cdots + {\bf X}_n ={\bf 0}) = 1.$$
\end{theorem}
{\bf Proof} Using Lemma 4.2 we can  choose  ${\bf X}\sim E_{p,n}({\bf 0}, {\bf \Sigma}\otimes {\bf \Phi}, {\Psi})$
  with all marginals  ${\bf X}_i$'s (the $i$th column of ${\bf X}, i=1,2,\cdots, n$) have the same $p$-elliptical distribution ${\bf E}_p({\bf 0},  {\bf \Sigma}, { \Psi})$, where
    ${\bf \Phi}\ge 0 $ is an  $n\times n$ matrix with diagonal elements are 1.
   Using Lemma 4.1  one finds that
${\bf X}_1 + \cdots + {\bf X}_n={\bf X}{\bf e}_n\sim  {\bf E}_n ({\bf 0},  {\bf\Sigma}\otimes ({\bf e_n^t\Phi e_n}),  \Psi)$.  Taking
  ${\bf \Phi}= (1-\rho)E_n +\rho {\bf e}_n{\bf e}_n^T$ with
$\rho=-\frac{1}{n-1}$, $E_n$ is $n\times n$  identity matrix.
 It follows that $$P({\bf X}_1 + \cdots + {\bf X}_n ={\bf 0}) = 1.$$ $\hfill\square$

\begin{corollary} If  $F\sim   {\bf E}_p({\boldsymbol \mu},    {\bf \Sigma},  \Psi)$,  where ${\bf \Sigma}\ge 0 $ is a  $p\times p$ matrix, $\Psi$  is a characteristic generator
for a    $p\times n$  matrix variate elliptically contoured distribution, then for any $n\ge 2$, there exist
 $p$-dimensional random vectors ${\bf X}_1,\cdots, {\bf X}_n$ identically distributed as $F$  such that
 $$P({\bf X}_1 + \cdots + {\bf X}_n =n{\boldsymbol\mu}) = 1.$$
\end{corollary}

\section{Conclusions and future work}
\setcounter{equation}{0}
 We  present three new proofs to a result due to  Wang, Peng and   Yang (2013) on JM of elliptical distributions
with the same characteristic generator.  We generalize this result to any distributions with finite second moments.
Moreover, we solve an  open problem proposed by  Wang (2015).  We also extent some results  to a  class of     multivariate   elliptically contoured distributions.
A full characterization of complete or joint mixability is still open.
In particular,  find necessary and sufficient conditions for complete mixability or joint mixability
of bounded  distributions or aymmetric distributions or multimodal distributions  there are still a lot of work to do.  Further open questions in this field   are collected in Wang (2015).

\noindent{\bf Acknowledgements.} \ We are grateful to Ruodu Wang  for insightful
comments  and valuable suggestions on the first manuscript.
This research   was supported by the National Natural Science Foundation of China (No. 12071251).

\end{document}